\def\ZZ         {{\bf Z}}
\def\RR         {{\bf R}}
\def\CC         {{\bf C}}
\def\QQ         {{\bf Q}}

\documentstyle[twoside,12pt]{article}
\newtheorem{prop}{Proposition}[section]
\newtheorem{dfn}[prop]{Definition}

\newtheorem{rem}[prop]{Remark}
\newtheorem{coro}[prop]{Corollary}

\title{String Cohomology of a Toroidal Singularity}

\author{
Lev A. Borisov
\\
\small Department of Mathematics,  Columbia University \\
\small 2990 Broadway, Mailcode 4432, New York, NY 10027, USA \\
\small e-mail: lborisov@math.columbia.edu}

\begin{document}

\date{}

\maketitle

\begin{abstract}
{We construct explicitly regular sequences in the semigroup ring
$R=\CC[K]$ of lattice points of the graded cone $K$. We conjecture that
the quotients of $R$ by these sequences describe locally string-theoretic
cohomology of a toroidal singularity associated to $K$. As a byproduct,
we give an elementary proof of the result of Hochster that semigroup 
rings of rational polyhedral cones are Cohen-Macaulay.
}
\end{abstract}

\section{Introduction}

String cohomology vector space of a variety $X$ with Gorenstein
toroidal singularities is a rather mysterious object. It is supposed
to be a chiral ring of no less mysterious $N=(2,2)$ superconformal
field theory constructed from $X$ and it has known graded dimension.
However, the space itself has not been identified so far in such
generality. The goal of this paper is to present a candidate
for the "contribution of a singular point" to this cohomology space.

The paper is organized as follows. Section 2 contains 
important preliminary results on the structure of lattice points
of the graded cone $K$. Section 3 uses these results 
to show that some explicitly written sequences of elements of
$R=\CC[K]$ are regular in $R$ and in $R$-module $R^{open}=\CC[K^{open}]$.
It also contains the proof of an analog of Poincar{\'e} duality.
It is worth mentioning that we give a short elementary proof of the theorem
of Hochster \cite{hochster}. Finally, the last section describes     
the relation of these results to Mirror Symmetry and string cohomology.

The author was inspired by recent preprints of Hosono \cite{hosono}
and Stienstra \cite{stienstra} who clarified the relationship between
the solutions of GKZ hypergeometric system and Mirror Symmetry.
The construction of this paper belongs to the A-side of Mirror Symmetry,
any B-side construction should involve solutions of GKZ systems.

One of the basic ideas of the argument has the flavor
of the theory of Gr{\"o}bner bases, which the author learned
from \cite{bayer.mumford}. It also appears that it involves 
the large complex structure limit, see for example \cite{morrison}.

Author would like to thank Dave Bayer and Sorin Popescu for helpful
remarks.

\section{Decomposition of Cone Lattice Points}

Let $N$ be a free abelian group of rank $r$. Let $K$ be a rational
polyhedral cone inside $N\otimes\RR$. We will assume that $K-K=N$ and
$K \cap (-K)=\{0\}$. We will also assume that the cone $K$ is graded,
that is there exists an element ${\rm deg}\in M={\rm Hom}(N,\ZZ)$ such that
the integer generators of all one-dimensional faces of the cone
$K$ have degree $1$. We will denote the interior of $K$ by $K^{open}$.

Another piece of data is a subset $\{e_i\},i=1,...,d$ of the set of lattice
points of degree $1$ that lie in $K$. The only condition on the subset is
that it includes the generators of all one-dimensional faces of $K$, that is 
$$\sum\RR_{\geq 0}e_i=K.$$ 
We also choose a maximum regular triangulation $T$ based on these points $e_i$
and denote by $\psi$ a strictly convex function on $K$ which
is linear on the simplices of triangulation $T$.

Our first goal is to construct a decomposition of the sets $K\cap N$
and $K^{open}\cap N$ into the disjoint union of sets $S_k$ of the 
form 
$$S_k=b_k+\sum_{i\in I_k}\ZZ_{\geq 0} e_i$$
where $I_k$ is a simplex of triangulation $T$ of maximum dimension $r$
and $b_k$ is a lattice point inside $\sum_{i\in I_k}\RR_{\geq 0}e_i$.


To carry out the construction for a given cone $K$ we fix a generic
vector $\xi\in N\otimes\RR$ that lies in $K^{open}$. For every 
$I\in T$ of maximum dimension, we consider the coordinates of $\xi$
in $I$, that is we look at $\beta_{I,i}$, such that
$$\xi=\sum_{i\in I}\beta_{I,i}e_i.$$
Because of the genericity of $\xi$, all $\beta$-s are non-zero.
We introduce the sets $B_{I,\xi}$ and $B_{I,-\xi}$ as follows
$$B_{I,\xi}=\{b\in I \cap N, ~{\rm such~that}~
b=\sum_{i\in I}\gamma_ie_i ~{\rm with}~
0<\gamma_i\leq 1 ~{\rm if}~\beta_{I,i}<0$$
$${\rm and}~ 0\leq\gamma_i< 1 ~{\rm if}~\beta_{I,i}>0\},$$
$$B_{I,-\xi}=\{b\in I \cap N, ~{\rm such~that}~
b=\sum_{i\in I}\gamma_ie_i ~{\rm with}~
0<\gamma_i\leq 1 ~{\rm if}~\beta_{I,i}>0$$
$${\rm and}~ 0\leq\gamma_i< 1 ~{\rm if}~\beta_{I,i}<0\}.$$

\begin{prop}
{\rm 
In the above notations the following statements hold.
  
(a) The set $K\cap N$ is the disjoint union of sets
$b+\sum_{i\in I}\ZZ_{\geq 0}e_i$ taken over all $I\in T$ of
maximum dimension and all $b\in B_{I,\xi}$.

(b) The set $K^{open}\cap N$ is the disjoint union of sets
$b+\sum_{i\in I}\ZZ_{\geq 0}e_i$ taken over all $I\in T$ of
maximum dimension and all $b\in B_{I,-\xi}$.
}
\label{disjoint}
\end{prop}

\noindent
{\em Proof.} (a) If $n\in K\cap N$, consider $n+\epsilon \xi$ for
very small $\epsilon>0$. It lies in $\sum_{i\in I}\RR_{\geq 0}e_i$
for some maximum simplex $I\in T$.
Therefore, we have  
$$n+\epsilon\xi = \sum_{i\in I}\alpha_i(\epsilon)e_i$$
and 
$$n = \sum_{i\in I}(\alpha_i(\epsilon)-\epsilon \beta_{I,i})e_i$$
where $\alpha_i(\epsilon)>0$.
Notice that $(\alpha_i(\epsilon)-\epsilon \beta_{I,i})$ are
independent of $\epsilon$. Therefore, they are always nonnegative.
Moreover, they are positive for such $i$ that $\beta_{I,i}<0$. 
This easily implies that $n\in b+\sum_{i\in I}\ZZ_{\geq 0}e_i$
for some $b\in B_{I,\xi}$.

Conversely, if $n\in b+\sum_{i\in I}\ZZ_{\geq 0}e_i$ with $b\in B_{I,\xi}$,
then for small $\epsilon>0$, the vector $n+\epsilon\xi$ lies
in $\sum_{i\in I}\RR_{>0}e_i$, which determines $I$ uniquely. 
Besides, there are clearly no intersections between 
$b_1+\sum_{i\in I}\ZZ_{\geq 0}e_i$ and $b_2+\sum_{i\in I}\ZZ_{\geq 0}e_i$
for different $b_1$ and $b_2$ from $B_{I,\xi}$.
The proof of (a) could now be finished by observation that
if $n+\epsilon \xi$ lies in some $\sum_{i\in I}\RR_{>0}e_i$ for small 
$\epsilon$, then $n$ lies in $K$.

(b) The proof is completely analogous. We use the fact that 
$n-\epsilon\xi$ lies in one of $\sum_{i\in I}\RR_{>0}e_i$
if and only if $n$ lies in $K^{open}$. \hfill$\Box$

\begin{coro}
{\rm 
Let us introduce polynomials 
$$S(t)=(1-t)^r\sum_{n\in K\cap N}t^{{\rm deg}(n)}~~{\rm and}~~ 
T(t)=(1-t)^r\sum_{n\in K^{open}\cap N}t^{{\rm deg}(n)}.$$
Then 
$$S(t)=\sum_{I,b\in B_{I,\xi}} t^{{\rm deg}(b)},~~
T(t)=\sum_{I,b\in B_{I,-\xi}} t^{{\rm deg}(b)}.$$
}
\label{SandT}
\end{coro}

\noindent
{\em Proof.} Follows immediately from the above proposition. \hfill$\Box$
\smallskip

Notice that the standard duality formula 
$$S(t)=t^rT(t^{-1})$$
follows immediately from this corollary together with the definitions of 
$B_{I,\xi}$ and $B_{I,-\xi}.$

In the next section we will use the following result.
Let us fix a lattice element $n\in K$. We look for all possible ways
of representing $n$ in the form 
$$n=b+\sum_{i=1}^d k_i e_i$$
where $k_i$ are non-negative integers and $b\in \cup B_{I,\xi}.$
The decomposition of $K\cap N$ above gives us one such representation
$$n=b_0+\sum_{i\in I_0}l_ie_i.$$
We claim that it has special properties
with respect to the convex function $\psi$.

\begin{prop}
{\rm
If $n=b+\sum_{i=1}^d k_i e_i$  then for small $\epsilon>0$
$$\psi(n+\epsilon\xi) \geq \psi(b+\epsilon\xi) + \sum_{i=1}^d k_i \psi(e_i)$$
and equality holds if and only if
$b=b_0,~k_i=l_i~{\rm for}~ i\in I_0,~k_i=0~{\rm for}~ i\notin I_0.$ 
}
\label{psiK}
\end{prop}
\noindent
{\em Proof.} The inequality is the basic property of the convex function
$\psi$. Equality holds if and only if there exists a maximum simplex $I$
such that the cone $\sum_{i\in I}\RR_{\geq 0} e_i$ contains $b+\epsilon\xi$
and all $e_i$ for which $k_i$ are non-zero. Therefore, 
$n+\epsilon\xi \in \sum_{i\in I}\RR_{\geq 0} e_i$ 
and the proof of Proposition \ref{disjoint}(a) shows that
$I=I_0$. Because of $b+\epsilon\xi \in \sum_{i\in I_0}\RR_{\geq 0} e_i$,
the lattice element $b$ lies in $b_1+\sum_{i\in I_0}\ZZ_{\geq 0}e_i$ for some
$b_1\in B_{I_0,\xi}$. Because the union of such sets is disjoint,
we have $b=b_1$. So $b\in B_{I,\xi}$ and therefore $b$ must equal $b_0$.
\hfill $\Box$
\smallskip

We will also need a similar statement for $K^{open}$.

\begin{prop}
{\rm Consider a lattice element $n$ in $K^{open}$.
If $n=b+\sum_{i=1}^d k_i e_i$ for some $b\in \cup B_{I,-\xi}$
then for small $\epsilon>0$
$$\psi(n-\epsilon\xi) \geq \psi(b-\epsilon\xi) + \sum_{i=1}^d k_i \psi(e_i)$$
and equality holds if and only if
$b=b_0,~k_i=l_i~{\rm for}~ i\in I_0,~k_i=0~{\rm for}~ i\notin I_0,$ 
where $n=b_0+\sum_{i\in I_0}l_ie_i$ is given by Proposition \ref{disjoint}(b).
}
\end{prop}

\noindent
{\em Proof.} The proof of this proposition is completely analogous
to the proof of the previous one. \hfill$\Box$

\bigskip

\section{Regular Sequences}

Let us fix a basis $m_1,...,m_r$ of the vector space 
$M\otimes \CC$ where $M={\rm Hom}(N,\ZZ)$.
We introduce the semigroup ring $R=\CC[K]$
and for every $n\in K$ we denote the corresponding element
in $R$ by $x^n$.
We also introduce $r$ elements of $R$ by the formula
$$Z_j = \sum_{i=1}^d <m_j,e_i> {\rm e}^{2\pi{\rm i} a_i}
x_i.$$
Here $a_i$ are some numbers assigned to the lattice elements
$e_i$ and the elements in $R$ that correspond to $e_i$
are denoted by $x_i$. 
These $Z_j$-s act on $R$ itself, and also on $R$-module
$R^{open}=\CC[K^{open}]$, which is an ideal in $R$.

The goal of this section is to show that for a generic choice of
$a_i$ the sequence $Z_1,Z_2,...,Z_r$ is regular on 
both $R$ and $R^{open}$. The following proposition is
crucial.

\begin{prop}
{\rm
Denote by $Z$ the ideal generated by $Z_1,...,Z_r$. Then
the following statements hold for generic $a_i$.

(a) Images of $x^b$ for $b\in \cup B_{I,\xi}$ generate
$R/ZR$ as $\CC$-vector space. 

(b) Images of $x^b$ for $b\in \cup B_{I,-\xi}$ generate
$R^{open}/ZR^{open}$ as $\CC$-vector space. 
}
\label{3.1}
\end{prop}

\noindent
{\em Proof.} (a) 
We introduce the ring $R_1=\CC[x_1,...,x_d]$ 
and consider $R$ and $R/ZR$ as $R_1$-modules.
Proposition \ref{disjoint} implies that 
these $R_1$-modules are generated by $x^b, ~b\in \cup B_{I,\xi}$.
Therefore, for each $q$ we have a surjective map
\begin{equation}
\oplus_{b\in \cup B_{I,\xi}} R_1[x^b] \to R/ZR \to 0 
\label{surj}
\end{equation}
of $R_1$-modules. 

The kernel of map (\ref{surj}) contains generators of two types.

\noindent $\bullet$ {\em Binomial relations.} 
Whenever we have an identity in the lattice $N$
$$
n=b_1+\sum_{i=1}^d k_{i1}e_i=b_2+\sum_{i=1}^d k_{i2}e_i
$$
we have a generator of the form 
$$
\prod_{i=1}^d x_i^{k_{i1}}[x^{b_1}]-\prod_{i=1}^d x_i^{k_{i2}}[x^{b_2}].
$$

\noindent $\bullet$ {\em Linear relations.}
We have generators $Z_jr_1[x^b]$ for $j=1,...,d,
~b\in \cup B_{I,\xi},~r_1\in R_1$.

It is enough to show that $\oplus\CC[x^b]$ maps surjectively
on the part of $R/ZR$ of degree less than some fixed big
number $D$. Really, it is enough to show that any element
of form $x_i[x^b]$ can be re-expressed as $\sum_b\alpha_b[x^b]$
modulo above relations, and degrees of $x^b$ are less
than $r$.

Let us pick a parameter $q$ and choose 
$${\rm e}^{2\pi{\rm i} a_i}=q^{\psi(e_i)}.$$
We will also make the following change of variables for each
non-zero $q$. We introduce
$$
(x_i)_{new}=q^{\psi(e_i)}x_i,~~[x^b]_{new}=q^{\psi(b+\epsilon\xi)}
[x^b]
$$
where $\epsilon$ is chosen to be small enough to fit in Proposition
\ref{psiK} for all $n$ of degree less than $D$. Then we rewrite the
generators of the kernel of map (\ref{surj}) in terms of new variables.

\noindent $\bullet$ {\em Binomial relations.}
Whenever we have an identity in the lattice $N$
$$
n=b_1+\sum_{i=1}^d k_{i1}e_i=b_2+\sum_{i=1}^d k_{i2}e_i
$$
we have a generator of the form 
$$
q^{\psi(b_1+\epsilon\xi)+\sum_i k_{i1}\psi(e_i)-
\psi(b_2+\epsilon\xi)-\sum_i k_{i2}\psi(e_i)}
\prod_{i=1}^d (x_i)_{new}^{k_{i1}}[x^{b_1}]_{new}
-\prod_{i=1}^d(x_i)_{new}^{k_{i2}}[x^{b_2}]_{new}.
$$

\noindent$\bullet$ {\em Linear relations.}
We have generators 
$$Z_jr_1[x^b]_{new}=
\sum_{i=1}^d <m_j,e_i>(x_i)_{new}r_1[x^b]_{new}.$$

Among the binomial relations, we will pick only the
ones where $n=b_1+\sum_{i=1}^d k_{i1}e_i$ is given by
the decomposition of Proposition \ref{disjoint}. Then,
by Proposition \ref{psiK}, the power of $q$ is
positive, unless $b_2=b_1,~k_{\cdot 2}=k_{\cdot 1}$.

Pick a basis of $\oplus_b (R_1[x^b]_{new})_{{\rm deg}<D}$
that consists of the products of monomials in $R_1$ and
$[x^b]_{new}$. For every $q$ we can introduce a matrix 
$A(q)$ which describes the map to 
$\oplus_b (R_1[x^b]_{new})_{{\rm deg}<D}$ from
the direct sum 
$$\oplus_b \CC[x^b]_{new}
~\oplus_{binomial} \CC[binomial]
~\oplus_{j,b,r_i} \CC Z_jr_i[x^b]_{new}
$$ 
where the direct sum is over the binomial relations that we have just
picked and $r_i$ are chosen to be monomials in $x_{new}$ of degree
less than $D$. 

To show that the vector space $\oplus_b \CC[x^b]_{new}$ surjects
onto $(R/ZR)_{{\rm deg}<D}$, it is enough to demonstrate that the matrix
$A(q)$ has full rank. Notice, that we have picked relations
in such a way that $A(q)$ has a limit $A(0)$ as  $q\to 0$.
Therefore, it will be enough to show that $A(0)$ has full
rank. 

The {\em binomial} relations become {\em monomial} in the
limit $q\to 0$ and hence the image of $A(0)$ contains
all basis elements of $\oplus_b (R_1[x^b]_{new})_{{\rm deg}<D}$ except,
perhaps, the elements of the form
$\prod_{i\in I} (x_i)_{new}^{k_i}[x^b]_{new}$
for $b\in B_{I,\xi}$.
However, if we use the {\em linear} relations, we 
can express $(x_i)_{new},i\in I$ in terms of other $(x_i)_{new}$, which 
shows that all the basis elements except for $[x^b]_{new}$ themselves
are in the image of $A(0)$. And since $[x^b]_{new}$ are also included
in the image of $A(0)$ by construction, we have the desired surjectivity 
of $A(0)$, which finishes the proof of (a).

The proof of (b) is completely analogous.
\hfill$\Box$
\smallskip

From now on we assume that $a_i$ are generic.
It is easy now to prove that $Z_1,...,Z_r$ form a regular
sequence on $R$ and $R^{open}$. We thus reprove for graded cones the
result of Hochster \cite{hochster} which states that $R$ is Cohen-Macaulay.

\begin{prop}
{\rm 
The sequence $Z_1,...,Z_r$ is regular on $R$ and $R^{open}$.
Thus $R^{open}$ is a Cohen-Macaulay module over the Cohen-Macaulay
ring $R$.
}
\label{regularR}
\end{prop}

\noindent
{\em Proof.} Let us show that $Z_1,...,Z_r$ is regular on $R$.
For every two power series $f(t)$ and $g(t)$ we say that $f(t)>g(t)$
if the first non-zero coefficient of $f(t)-g(t)$ is positive.

For each $k=0,...,r$ we denote 
$$f_k(t)=\sum_{l\geq0}t^l {\rm dim}_{\CC} (R/(Z_1,...,Z_k)R)_{{\rm deg}=l}
~.$$ 
The exact sequence 
$$R/(Z_1,...,Z_k)R \to R/(Z_1,...,Z_k)R \to R/(Z_1,...,Z_{k+1})R
\to 0$$
implies that $$f_{k+1}(t)\geq f_k(t).$$

On the other hand, the fact that $\oplus_{b\in \cup B_{I,\xi}} \CC[x^b]$
surjects onto $R/ZR$ implies that the power series (in fact, it is a
polynomial) $f_r(t)$ is less or equal to $\sum_{I,b\in B_{I,\xi}} t^{{\rm deg}(b)}$
which is equal to $(1-t)^r f_0(t)$ by Corollary \ref{SandT}.
Therefore, all intermediate inequalities are equalities, which
shows that the above sequences are exact on the left. 

The same argument works for $R^{open}$. \hfill$\Box$

\begin{rem}
{\rm Theorem of Hochster could be proved in full generality using our methods.
Really, for any cone $K$ we can pick points $e_i$ on one-dimensional faces that
lie in the same hyperplane ${\rm deg}=1$ for some ${\rm deg}\in M\otimes\QQ$. 
Then the only difference is that ${\rm deg}(n)$ is allowed to take values
in ${1\over l}\ZZ$ for some $l$, which also requires the use
of fractional powers of $t$. However, this does not present any problems,
because the integrality of ${\rm deg}(n)$ was never used.
}
\end{rem}

\begin{coro}
{\rm
Surjective maps of Proposition \ref{3.1} are isomorphisms.
}
\end{coro}

\noindent
{\em Proof.} It follows from the proof of Proposition \ref{regularR}
that graded dimensions of these spaces are the same, so surjectivity
implies bijectivity. \hfill$\Box$
\smallskip

\begin{rem}
{\rm
Regularity of the sequence $Z$ was used in a special case without proof 
in the paper \cite{GKZ}. In the later correction note \cite{GKZcorr}
the result is stated explicitly, but the proof is inadequate.
}
\end{rem}
\smallskip

Because of the duality $S(t)=t^rT(t^{-1})$, we have 
${\rm dim}_\CC (R^{open}/ZR^{open})_{{\rm deg}=r}=1$. We denote by
$\varphi$ a surjective map $R^{open}/ZR^{open}\to \CC$ which sends
$(R^{open}/ZR^{open})_{{\rm deg}<r}$ to zero. Then we have a pairing
$$(R/ZR)\otimes_{\,\CC}(R^{open}/ZR^{open}) \to \CC$$
which maps $x\otimes y$ to $\varphi(xy)$.

\begin{prop}
{\rm 
({\em Poincar{\'e} Duality}\/)
The pairing $$(R/ZR)\otimes_{\,\CC}(R^{open}/ZR^{open}) \to \CC$$
is non-degenerate.
}
\label{Poincare}
\end{prop}

\noindent
{\em Proof.}
We need to show that for every element $x\in R^{open}/ZR^{open}$
the principal $R$-submodule it generates inside $R^{open}/ZR^{open}$
is non-zero at degree $r$. Let us pick a homogeneous $x$ whose 
principal submodule is zero in degree $r$, which has the highest
degree (less than $r$) among all $x$ with this property.
Denote by $R_{>0}$ the maximum ideal in $R$. For every homogeneous 
$y\in R_{>0}$ the principal submodule of $xy$ is zero in degree $r$,
but $xy$ has a higher degree, 
so it must be zero. This implies that there is a non-trivial homomorphism 
from $\CC=R/R_{>0}$ to $R^{open}/ZR^{open}$, which maps $1$ to $x$.
Since the top element certainly provides us 
with a homomorphism $\CC\to R^{open}/ZR^{open}$, it suffices to
show that
$${\rm Hom}^R(\CC,R^{open}/ZR^{open})\cong\CC.$$

Now we use the well-known result (see, for example, \cite{danilov})
that $R^{open}$ is the canonical module for $R$. Hence
$${\rm Ext}^R_i(\CC,R^{open})\cong 0, ~i\neq r,
~{\rm Ext}^R_r(\CC,R^{open})\cong\CC,$$
which is a standard property of canonical modules,
see \cite{bruns.herzog}.
Now it can be easily deduced from the Koszul complex
associated to $Z$ and $R^{open}$ that 
$${\rm Hom}^R(\CC,R^{open}/ZR^{open})\cong{\rm Ext}^R_r(\CC,R^{open}) \cong \CC,$$
which completes the proof. \hfill$\Box$

\bigskip
\section{Relation to Mirror Symmetry and Other Comments}

Now it is time to explain the title of the paper. String-theoretic
cohomology of a variety $X$ with toroidal Gorenstein singularities is supposed
to be the chiral ring of the corresponding $N=(2,2)$ superconformal
field theory. It is still not clear how to construct such a theory
for varieties with most general singularities of this type.
However, graded dimension of string-theoretic cohomology
vector spaces was suggested by Batyrev and Dais in their paper 
\cite{batyrev.dais}. It was later verified in \cite{bat.bor}
that this definition of string-theoretic Hodge numbers is compatible 
with mirror duality of Calabi-Yau hypersurfaces and complete intersections
in Gorenstein toric Fano varieties.

\begin{dfn}
{\rm \cite{batyrev.dais}
 Let $X = \cup_{i \in I} X_i$ be a stratified algebraic variety 
over $\CC$ with at most Gorenstein toroidal singularities 
such that for any $ i \in I$ the singularities of $X$ along 
the stratum $X_i$ of codimension $k_i$ are defined by a $k_i$-dimensional
finite rational polyhedral cone $K_i$; that is  $X$
is locally isomorphic to
$$\CC^{{\rm dim}(X)-k_i} \times U_{K_i}$$
at each  point $x \in X_i$ where $U_{K_i}$ is a
$k_i$-dimensional  affine toric variety which is associated
with the cone $K_i$ (see \cite{danilov}).
Batyrev and Dais have introduced the polynomial
$$
E_{\rm st}(X;u,v) = \sum_{i \in I} E(X_i;u,v) \cdot
S_{K_i}(uv)
$$
where $E(X_i;u,v)$ are E-polynomials of Danilov and Khovanski\^i,
see \cite{danilov.khov}.
It is called the {\em string-theoretic E-polynomial of $X$.} If we write
$E_{\rm st}(X; u,v)$ in form
$$
E_{\rm st}(X;u,v) = \sum_{p,q} a_{p,q} u^{p}v^{q},
$$
then the numbers $h^{p,q}_{\rm st}(X) = (-1)^{p+q}a_{p,q}$
are called the {\em string-theoretic Hodge numbers of $X$.}
}

\end{dfn}

Thus it seems that the "local" description of string cohomology should
be provided by a vector space that has graded dimension $S(t)$.
This is precisely what we have achieved in the previous sections. We now
suggest that the vector space $R/ZR$ should be thought of as the
local contribution of a toroidal singularity to string cohomology of
the variety $X$. There also seems to be a nice notion of local string 
homology, which is provided by $R^{open}/ZR^{open}$. The existence of 
natural pairing that satisfies Poincar\'e duality provides a further 
justification of our terminology.

It is important to remark that we had to choose some numbers 
$a_i$ to facilitate the construction, and that the resulting spaces
do depend on this choice. This is due to the fact that the superconformal
field theory in question depends on both complex and K{\"a}hler parameters
when $X$ is smooth. So it is natural to suggest that $a_i$ play the role
of K{\"a}hler parameters here. This is best illustrated by the example
of hypersurfaces in mirror dual Gorenstein toric Fano varieties.
To construct an $N=(2,2)$ theory, we really need a pair of such hypersurfaces,
and it roughly amounts to choosing coefficients $a_i$ for points
on both dual reflexive polyhedra.

In general, if we choose consistently the numbers $a_i$ for integer points
of degree one in cones $K_i$, we can define a (non-coherent) 
sheaf ${\cal A}$ of ${\CC}$-algebras over $X$ as follows. 
Whenever the closure of one strata $X_i$ contains
another strata $X_j$ there is  a surjective map
$$(R/ZR)_{K_j} \to (R/ZR)_{K_i}.$$
We define the germ ${\cal A}_x$ over a point $x\in X_i$ to be 
$(R/ZR)_{K_i}$ and sections are, by definition, constant on each strata 
and compatible with the above maps. 
The sheaf ${\cal A}$ is naturally graded. One can hope 
to somehow use the cohomology of this sheaf to define 
string-theoretic cohomology vector spaces.
The cohomology of ${\cal A}$ appears to be a viable
candidate for string cohomology of $X$ when it is an orbifold. 
Unfortunally, in general the Hodge structure of $H^*({\cal A})$ 
is not pure, which is certainly one of
the properties to be expected of string cohomology.

As a side remark, the construction of this paper 
works not only for $R$ and $R^{open}$
but for some other ideals of $R$ that are associated with the choice of 
$\xi$ which is neither in $K$ nor in $-K$. It would be interesting to see 
whether these ideals have any additional nice properties.

The construction here has the flavor of the A-side of Mirror Symmetry.
It would be very interesting to see the B-side construction, which is
presumably more useful. It is very possible, that the spaces constructed
here could be mapped to the spaces of solutions of GKZ hypergeometric
systems and their duals. One may also try to define a flat connection on
the vector bundle with fibers $R/ZR$ over the space of parameters $a_i$.
This is the direction of further research that the author plans to pursue.

\end{document}